\documentclass[a4paper,10pt]{article}
\usepackage[all]{xy}
\usepackage{amsmath}\usepackage[latin1]{inputenc}
\usepackage{epsfig}\usepackage{amsthm}\usepackage{amssymb}\usepackage{epsfig}\usepackage{amscd}
\usepackage{graphicx,subfig,xcolor}
\usepackage{fancyhdr}
\theoremstyle{plain}
\newtheorem{thm}{\bf Theorem}[section]
\newtheorem{prop}[thm]{\bf Proposition}
\newtheorem{lem}[thm]{\bf Lemma}
\newtheorem{cor}[thm]{\bf Corollary}
\newtheorem{de}[thm]{\bf Definition}
\newtheorem{ex}[thm]{\bf Example}
\newtheorem{no}[thm]{\bf Notation}

\title{Prequantization of logsymplectic structures.}
\author{J.DONGHO}

\begin{document}

\maketitle

\begin{abstract}
In this paper, we study quatization condition of logsymplectic structure using integrality of such structue
on the complement of associated divisor $D$.
\end{abstract}

\section{Introduction}
The notion of logsymplectic structure taking their origin from particular meromorphic
forms having at most simple poles along certain divisor $D$ of a giving complex manifold $X.$
 Such forms are amply studies in \cite{KS}. Using the notion Lie-Rinehart algebra, we give the
 algebraic generalization of such notion. Giving an $K$ algebra $\mathcal{A},$ and $I$ a non trivial
ideal of $\mathcal{A},$ we consider the $\mathcal{A}$-submodule $Der_k(\log I)$ of $Der_k(\mathcal{A})$
constituted by $\delta\in Der_k(\mathcal{A})$
such that $\delta(I)\subset I.$ It is usually called the $\mathcal{A}$-module of derivation logarithmic
 along $I.$ An element $\delta$ is called principal if for all $u\in I, \delta(u)\in u\mathcal{A}.$
 We denoted $\widehat{Der_k(\log I)}$ the subset of principal elements of $Der_k(\log I).$ It is a submodule of $Der_k(\log I)$ closed under $\mathcal{A}-$ module and $k-$ Lie structure of $Der_k(\mathcal{A}).$ The dual of $\widehat{Der_k(\log I)}$ is the $\mathcal{A}-$module $\Omega_k(\log I)$ generated by $\Omega$ and the set $\{\dfrac{du}{u}; u\in I-\{0\}\}$ The inclusion of $\widehat{Der_k(\log I)}$ in $Der_k(\log I)$ induce an algebra Lie-Rinehart structure on $\widehat{Der_k(\log I)}.$ We can then talk about the notion of Lie-Rinehart Poisson structure on $\widehat{Der_k(\log I)}$; it is a 2-form $\eta:\widehat{Der_k(\log I)}\otimes \widehat{Der_k(\log I)}\rightarrow \mathcal{A}$ such that $d_\rho\eta=0;$ where $d$ is the Chevaley lie-Rinehart differential of $\widehat{Der_k(\log I)}.$ A lie-Rinehart Poisson structure $\eta$ is called logsymplectic if it is non degenerated; in that case, $\mathcal{A}$ is called a logsymplectic algebra.

When $I$ designed the ideal of a giving divisor of some complex manifold $X,$ a
Lie-Rinehart-symplectic structure on $Der_\mathbb{C}(\log I)$ correspond to the notion of
 logsymplectic structure using by R. Goto in \cite{RG}.
In this particular case $Der_\mathbb{C}(\log I)$ is equal to $Der_X(\log D)$.

Recalled that giving a symplectic Lie-Rinehart algebra $(L, \rho, \eta)$ and $a\in\mathcal{A},$
there is a unique element $\delta_a\in L$ such that
\begin{center}
 $i_{\delta_a}\eta=d_\rho a$
\end{center}
Let $\mathcal{L}_\delta=[i_\delta, d_\rho]$; where $[,]$
denote the commutator and $\delta\in L,\mathcal{L}_{\delta_a}\eta=0$.
 From unicity of $\delta_a$, the following bracket is well defined,$\{a,b\}=\eta(\delta_a,\delta_b).$
 which is Poisson structure induced by Lie-Rinehart symplectic structure $\eta.$ When the Lie-Rinehart symplectic algebra is $(Der_k(\log I),\eta),$ the associated Poisson structure is called logsymplectic Poisson structure.

Logsymplectic Poisson structure represented a particular case of degenerated Poisson structure since
 there are symplectic on the complement of the associated divisor.

The goal of this paper is to study integral condition of such structure. \\
According to Izu Vaisman in 
a Poisson manifold $(X, P$ is prequantization if and only if there exist vector field $X$ and 2-form
 $\omega$ who represented an integrable class such that $P+[X,P]_{SC}=P^{\sharp}(\omega);$ where $[-,-]_{SC}$
is well known Schouten bracket.  We can remark that Vaisman condition presented 2 main difficulties:
solve a partial differential equation in 2 variables $P+[X,P]_{SC}+ Y=0$ and compute the Poisson cohomology class of inverse of $Y$ for a giving solution $(X,Y).$  Follow B.Kostant and Soureau, Vaisman condition is equivalent to integral condition of $\eta$ in the symplectic case. The integral condition only involve the De Rahm cohomology class of $\eta.$ It will be very useful to know how we can change the Vaisman condition by someone filling integral condition; when we have singular Poisson structure.

In this paper, we study the case of singular Poisson structure;
 a logsymplectic Poisson structure recalled above.\\ Of it, we change the De Rham cohomology by logarithmic De Rham cohomology\\
Due to the fact that logsymplectic structure are symplectic on the complement of the divisor $D$,
we apply integral condition on the complement of $D$ and extend the corresponding prequantum line bundle to $X.$
\section{Preliminaries}
We begin this section by introducing the notion of
Lie-Rinehart-Poisson-Logsymplectic algebra which allow us to give
the algebraic definition of logsymplectic structure. This notion
breathe in \cite{JH}. It is particular case of
Lie-Rinehart-Poisson-symplectic algebra fully study in algebra and
Poisson geometry.
\subsection{On Lie-Rinehart Poisson-logsymplectic algebra.}
 In what follow, $\mathcal{A}$ is associative, commutative and unitary $k$-algebra on
a field $k$ such that $car(k)=0$. Let $I$ be a nonzero ideal of
$\mathcal{A}.$ $Der_k(\mathcal{A})$ is the $\mathcal{A}$-module of
derivations on $\mathcal{A}$ and $Der_k(\log I)$ is his submodule
constitute by logarithmic principal derivations along $I.$ We recall
that a Lie-Rinehart algebra is a pair $(L, \rho)$; where $L$ is Lie
algebra and $\rho:L\rightarrow Der_k(\mathcal{A})$ is Lie algebras
homomorphism satisfy the following equality.
\begin{equation}\label{Eq1}
    [l_1,al_2]=\rho(l_1)(a).l_2 + a[l_1,l_2]
\end{equation}
It follow from this definition that $Der_k(\mathcal{A})$ endowed
with identity is Lie-Rinehart algebra. In other hand, we can easily
prove that the inclusion map of $Der_k(\log I)$ endowed it to the
structure of Lie-Rinehart algebra. We can then consider on
$\mathcal{L}_{alt}(Der_k(\log
I),\mathcal{A})=\underset{{\mathbb{Z}}}\bigoplus{\mathcal{L}_{alt}}^p(Der_k(\log
I),\mathcal{A})$ the following differential:
\begin{equation}\label{Eq2}
\begin{array}{ccc}
    (d^{log}f)(\delta_0,...,\delta_p)=\underset{{i=1}}{\overset{p}\sum}(-1)^{i+1}\delta_i.f(\delta_1,...,\hat{\delta_i},...,\delta_p)
    \\+ \underset{{i<j}}{\sum}(-1)^{i+1}f([\delta_i,\delta_j],\delta_1,...,\hat{\delta_i,}...,\hat{\delta_j},...,\delta_p)
\end{array}
\end{equation}
As in general case, we have ${d^{\log}}^2=0;$ then corresponding
cohomology is called logarithmic De Rham cohomology. This cohomology
is fully study in the framework of algebraic and complex geometry.
For example,if $I$ denote the definition ideal of reduced divisor
$Y$ of a complex manifold $X$, $Der_{\mathbb{C}}(\log I)$ correspond
at each point $x$ to the $\mathcal{O}_{X,x}$-module of logarithmic
vector field introduce in \cite{KS}which is denoted $Der_{X,x}(\log
Y)$ and ${\mathcal{L}_{alt}}^p(Der_k(\log
I),\mathcal{A})={\Omega^p}_X(\log Y).$\\ Since $Der_k(\log I)$ is
Lie-Rinehart algebra, we can define on the notion of
Lie-Rinehart-Poisson structure. It is 2-cocycle of $d^{log}.$ Now,
we can give the definition of Lie-Rinehart-Poisson-logsymplectic
structure.
\begin{de}
{\sf A lie-Rinehart-Poisson-logsymplectic structure is a pair\\
$(Der_k(\log I),\mu)$ where $\mu$ is a non degenerated 2-cocycle of
$d^{\log}$}
\end{de}
As in smooth case, Lie-Rinehart-Poisson-logsymplectic structure
induce a Poisson structure on $\mathcal{A}.$ Indeed, since the
2-cocycle $\mu$ is non degenerated, its contraction by logarithmic
derivation induce an isomorphism $i$ of $\mathcal{A}$-modules from
$Der_k(\log I)$ to its algebraic dual $Der_k(\log I)^*$. So,
according to Hochschild-Konstant-Rosenberg \cite{HKR}, since
$Der_k(\mathcal{A})$ is submodule of $Der_k(\log I)$, for all
$a\in\mathcal{A}, d^{\log}(a)\in \Omega_k$ and
$\Omega_k=Der_k(\mathcal{A})$ if and only if $\mathcal{A}$ is
regular affine algebra on perfect field $k.$ In this case, we are
sure that $d^{\log}(a)\in Der_k(\log I)^*$ for all
$a\in\mathcal{A}.$ Then there is an unique $\delta_a\in Der_k(\log
I)$ such that:
\begin{equation}\label{Eq3}
    i_{(\delta_a)}\mu=d^{\log}(a).
\end{equation}
$a$ is called Hamiltonian of Hamiltonian derivation $\delta_a.$\\
For all $a, b\in\mathcal{A}$,  consider:
\begin{equation}\label{Eq4}
    \left\{ a,b\right\}=-\mu(\delta_a,\delta_b)
\end{equation}
We have have the following proposition.
\begin{prop}
Let $\mathcal{A}$ be a regular affine algebra on perfect field $k$
endowed with a Lie-Rinehart-Poisson-logsymplectic structure $\mu$.
 $\left\{ a,b\right\}=-\mu(\delta_a,\delta_b)$ is a well defined
 logarithmic Poisson structure on $\mathcal{A}.$
\end{prop}
The hypotheses of the above proposition is satisfy when
$\mathcal{A}$ is the sheaf of holomorphic functions on logsymplectic
manifold $X.$ Indeed, It is proven it \cite{KS} that the pair
$\{Der_X(\log Y), \Omega_X(\log Y)\}$ is reflexive. Then for all
$f\in\mathcal{O}_X, d^{\log}(f)\in \Omega_X(\log Y)=Der_X(\log
Y)^*.$ An then each $f\in\mathcal{O}_X$ is associated to an unique
$\delta_f\in Der_X(\log Y).$ Therefore, the above bracket is well
defined.\\ More generally, if the pair $\{Der_k(\log I),
\Omega_k(\log I)\}$ is reflexive, the above proposition is true and
besides, for all $u\in I-{0}, \dfrac{d^{\log}u}{u}\in \Omega_k(\log
I)=Der_k(\log I)^*.$ Therefore, there exist an unique
$\tilde{\delta}_u\in Der_k(\log I)$ such that
$i_{\tilde{\delta}_u}\mu=\dfrac{d^{\log}u}{u}.$ In this case, since
$I$ is subset of $\mathcal{A}$, for all $u\in I$ there exist
$\delta_u$ such that $i_{\delta_u}\mu=d^{\log}u.$ It is easy to
prove that $\delta_u=u\tilde{\delta}_u.$ We can then consider the
following bracket:
\begin{equation}\label{Eq5}
\left\{ a, b\right\}_{sing}:=\left\{
\begin{array}{ccccc}
\dfrac{1}{uv}\{u,v\}&\texttt{if}& a=u, b=v\in I-0\\\\
\dfrac{1}{u}\{u, b\} &\texttt{if}& a=u\in I-0, b\in\mathcal{A}-I\\\\
\{a, b\} &\texttt{if}& a,b\in\mathcal{A}-I
\end{array}\right.
\end{equation}
We have the following proposition.
\begin{prop}
If the pair $\{Der_k(\log I), \Omega_k(\log I)\}$ is reflexive then
each Lie-Rinehart-Poisson-Logsymplectic structure on $\mathcal{A}$
induce two Lie structure $\{-,-\}$ and $\{-,-\}_{sing}$ such that
for all $u,v\in I-{0}$,
\begin{itemize}
    \item [i)] $i_{({\delta_{\{u,v\}}-uv\delta_{\{u,v\}_{sing}}})}\mu=\{u,v\}(\dfrac{d^{\log}u}{u}+\dfrac{d^{\log}v}{v})$
    \item [ii)] $\{uv,a\}_{sing}=\{u+v,a\}_{sing}; \forall a\in\mathcal{A}-I.$
    \item [iii)] $\{a,b\}=\delta_a(b)$
    \item [iv)] $[\delta_a,\delta_b]=\delta_{\{a,b\}}$
    \item [v)] $\delta_{\{u,v\}}= uv[\tilde{\delta_u},\tilde{\delta_v}]+\{u,v\}(\tilde{\delta_v}+\tilde{\delta_u})$
\end{itemize}
\end{prop}
\begin{proof}It straightforward \end{proof} Let
 $\mathcal{H}^{\log}(\mathcal{A},I):=\{\delta\in Der_k(\log I); \exists a\in \mathcal{A};
 \delta=\delta_a\}.$ According to above proposition,
 $\mathcal{H}^{\log}(\mathcal{A},I)$  is sub-Lie-algebra of $Der_k(\log
 I)$ and we have the following exact sequence of Lie algebras.
 \begin{equation}\label{Eq6}
\xymatrix{O\ar[r]&k\ar[r]&(\mathcal{A},\{-,-\})\ar[r]^{\delta}&\mathcal{H}^{\log}(\mathcal{A},I)\ar[r]&0}
 \end{equation}
According to Dirac principe of quantization, $(\mathcal{A},\{-,-\})$
is quantization if there exist an representation
$(\mathcal{H},\varphi)$ where $\mathcal{H}$ is Lie-Rinehart
extension of $\mathcal{A}$ along $Der_k(\log I)$ such that the
following diagram commute.
\begin{equation}\label{Eq7}
\xymatrix{0\ar[r]&\mathcal{A}\ar[r]^f&\mathcal{H}\ar[r]^g&Der_k(\log I)\ar[r]&0\\
         0\ar[r]&k\ar[r]\ar[u]&(\mathcal{A},\{-,-\})\ar[u]^{\varphi}\ar[r]&
         \mathcal{H}^{\log}(\mathcal{A},I)\ar[u]\ar[r]&0}\end{equation}
The morphism $\varphi$ is called quantization formula and it
satisfy:
\begin{equation}\label{Eq8}
\varphi (as)=\nabla_{v(a)}s+2i\pi as
\end{equation}

$a\in \mathcal{A}$, $\nabla$ is a section of
$g$; and  $v(a)=\{a, -\}.$\\
In general, $\nabla$ is only an $\mathcal{A}$-module homomorphism. Obstruction to become Lie-morphism is
measured by cohomology class of an 2-cocycle $K_{\nabla}$;
it is usually called curvature of $\nabla$ on $\mathcal{H}.$\\ When $(\mathcal{H},\varphi)$ exist and $\varphi$
satisfy \ref{Eq8} the triplet $(\mathcal{H},\nabla, K_{\nabla})$ is called prequantum representation of $\mathcal{A}.$
 The following paragraph is devoted to the construction of
$\mathcal{H}$ when $\mu:=\omega$ is a logsymplectic structure on a
even dimensional complex manifold $X$ with reduced divisor $D.$

\subsection{Logarithmic differential operator.}
In this paragraph, $(X,\omega, D)$ is logsymplectic manifold.
Our stumbles is to construct a
solution $(\mathcal{H}, \varphi)$ to the above problem when $\varphi$ is given by (\ref{Eq8}).\\
Our construction of $\mathcal{H}$ is motivated by the notion of logarithmic
connection which is sufficiently studied in Complex and Algebraic geometry. We will also denote
 $\mathcal{E}$ a locally free $\mathcal{O}_X$-module of rank
$1$ and $D=\{h=0\}$ a divisor of $X.$
\subsubsection{Logarithmic connection.}
The notion of logarithmic connection is original in the work of P. Deligne when he formulated and proved the theorem
establishing a Riemann-Hilbert correspondence between monodromy groups and Fuchsian
systems of integrable partial equations or flat connections on complex manifolds.
He also gave a treatment of the theorem of Griffiths which states that the Gauss-Manin or Picard-Fuchs systems of
of differential equations are regular singular.
\begin{de}{\sf  Let $\mathcal{M}$ be $\mathcal{O}_X$-module. A connection
  on $\mathcal{M}$
with logarithmic poles along $D$ is a $\mathbb{C}$-linear
homomorphism $\nabla$ \[\nabla: \mathcal{M}\rightarrow
\Omega^1_X(\log D)\otimes\mathcal{M}\] that satisfies Leibniz's
identity:
\begin{equation}\label{Eq9}
\nabla(fm)=df.m +f\nabla (m)\end{equation}
 where $d$ is the exterior
derivative over $\mathcal{O}_X$. }
\end{de}
This is equivalent to a linear map $\triangle:Der_X(\log D)\rightarrow End(\mathcal{E})$ satisfy the following
\begin{equation}\label{Eq10}
    \triangle_X(fs)=f\triangle_X s + X(f)s
\end{equation}

If  $\nabla$ is logarithmic connection $K_{\nabla}$
 will denoted its curvature and the pair $\left( \mathcal{M},\nabla\right) $
 will refer to logarithmic connection on a locally free $\mathcal{O}_X$-module of rank
$1$ $\mathcal{M}$.
\begin{lem}
{\sf If $\left( \mathcal{M},\nabla\right) $ is rank one connection
on $X$, Then for all closed 1-form $\tau\in H^0\left( X,
\Omega^1_X(\log D)\right) $ $\left( \mathcal M,\nabla + \tau\otimes
id\right) $ is a connection with curvature form $K=K_{\nabla}$
 }
\end{lem}
\begin{proof}
Suppose that $\nabla$ is define by
$$\nabla (s)=\sigma\otimes s$$ for a giving nowhere vanish section $s$
\begin{displaymath}
\begin{array}{lll}
(\nabla +\tau\otimes id)(s) & = & \nabla(s) +\tau\otimes s \\
 & = & \sigma\otimes s +\tau\otimes s \\
 & = & (\sigma +\tau)\otimes s
\end{array}
\end{displaymath}
And \begin{displaymath}
\begin{array}{lll}
(\nabla +\tau\otimes id)(fs) & = & \nabla (fs) +\tau\otimes id(fs) \\
 & = & df\otimes s +f\sigma\otimes s +f\tau\otimes s \\
 & = & df\otimes s + f(\nabla +\tau\otimes id)s
    \end{array}
    \end{displaymath}

\end{proof}

Let $s_0\in
H^0(U_i,\mathcal{M})$ such that $0\notin s_0(U_i).$ \\ There exist
$\sigma\in H^0(U_i, \Omega^1_X(\log D))$ such that $\nabla
s_0=\sigma\otimes s_0.$ Then $K_\nabla=d\sigma.$

\noindent Let $p$ be a point of $D$ and $(z^i_\lambda)$ logarithmic
coordinate system along $D$ at $p.$
\begin{equation}\label{E5}
\sigma=\underset{i=1}{\overset{r}\sum}a_i\dfrac{dz^i_\lambda}{{z^i}_{\lambda}}
 +\underset{i=r+1}{\overset{2n}\sum}a_id{z^i}_{\lambda}
\end{equation}
 where $a_i\in H^0(X,\mathcal{O}_X).$
 It follow that
\begin{lem}\label{L1}
Let $D$ be a normal crossing divisor and $\alpha\in H^0(X,\Omega_X(\log
D))$. If $d\alpha=0$ then the residue of $\alpha$ is constant on any
component of singular locus of $D.$ Any such form with at least one
nonzero residue admits representation
\begin{equation}\label{E6}
    \alpha=\underset{j=1}{\overset{r}\sum}\alpha_i\dfrac{df_j}{f_j},\quad\quad\quad
    \alpha_1,...,\alpha_r\in \mathbb{C}
\end{equation}
\end{lem}
\begin{proof}
Let $p$ be a point of $D$ and $U_\lambda$ an open coordinate
neighborhood of $p$. We have:\
$\alpha=Res(\alpha)\dfrac{dh_p}{h_p}+\alpha_{reg}$ Where
$Res(\alpha)$ is the residue of $\alpha.$ $d\alpha=0$ imply
$d(Res\alpha)=O$, (since $Res$ commute with $d$). However, From
Theorem 2.9 in \cite{KS}, $Res(\Omega^1_X(\log D))=\mathcal{O}_X$,
 $dRes(\alpha)=0$ imply $Res\alpha\in\mathbb{C}.$
\end{proof}
  \begin{prop}{\sf $K_\nabla =0$ if and only if $\sigma=\underset{i=1}{\overset{r}\sum}a_i
  \dfrac{dz^i_\lambda}{{z^i}_{\lambda}}$ with $a_i\in\mathbb C$}
 \end{prop}

\begin{proof}
$K_\nabla=d\sigma$; where $\sigma$ is the connection one form of
$\nabla.$ The result is then a consequence of Lemma \ref{L1}
\end{proof}

\begin{de}{\sf
Let $(\mathcal{M}, \nabla)$ and $(\mathcal{N},\delta)$ be two
connections. An homomorphism from $(\mathcal{M}, \nabla)$ to
$(\mathcal{N},\delta)$ is a sheaf homomorphism
$\varphi:\mathcal{M}\rightarrow\mathcal{N}$ such that the following
diagram commute.
\[\xymatrix{\mathcal{M}\ar[r]^{\varphi}\ar[d]_{\nabla}&\mathcal{N}\ar[d]^{\delta}\\
\Omega^1_X(\log
D)\otimes_{\mathcal{O}_X}\mathcal{M}\ar[r]&\Omega^1_X(\log
D)\otimes_{\mathcal{O}_X}\mathcal{N}}\]

 }\end{de}
\begin{de}{\sf Let $(\mathcal{M},\nabla)$ be a connection on
$X^*$. A meromorphic prolongation of $(\mathcal{M},\nabla)$ is a
meromorphic connection $(\bar{\mathcal{M}},\bar{\nabla})$ on $X$
such that the restriction is an isomorphism.}\end{de}
\subsubsection{Module of logarithmic differential operator.}
 Let $\mathcal{A}$ be a commutative ring. For a pair of $\mathcal{A}$-modules $M,N$ we define module $Diff_\mathcal{A}^k(M,N)$ inductively by putting
\begin{enumerate}
 \item $Diff_\mathcal{A}^0(M,N)=Hom_\mathcal{A}(M,N)$
\item $Diff_\mathcal{A}^k(M,N)=\\
\left\lbrace \texttt{additive maps}
 \,\, u:M\rightarrow N \,\,\texttt{s.t}\, \forall a\in\mathcal{A}\, au-ua\in Diff_\mathcal{A}^{k-1}(M,N) \right\rbrace $
\end{enumerate} Elements of $Diff_\mathcal{A}^k(M,N)$ are called $k$-order differential operator from $M$ to $N.$ We note $Diff_\mathcal{A}^k(M)$ for $Diff_\mathcal{A}^k(M,M).$

Replacing $M$ by $\mathcal{E},$ the above definition became;
 \begin{de}\label{Def 2.4}{ \sf A r-order differential operator on $\mathcal{E}$ is a
 $\mathbb{C}$-linear map $\varphi:\mathcal{E}\rightarrow\mathcal{E}$
 such that $s\mapsto\varphi(fs)-f\varphi (s)$ is an (r-1)-order differential operators on
 $\mathcal{E})$; for all $f\in \mathcal{O}_X$}
 \end{de}
 In the previous paragraph, we see that each logarithmic connection induce a morphism $\triangle$ from $Der_X(\log D)$
 such that for all $f\in\mathcal{O}_X,X\in Der_X(\log D), \triangle_X(fs)-f\triangle_X(s)=X(f)s.$
 It follow that $\triangle_X\in End(\mathcal{A})$ and the map $\triangle_X^f: s\mapsto \triangle_X(fs)-f\triangle_X(s)$
 is zero order differential operator on $\mathcal{E}$ and $\triangle_X^f(h)\in h\mathcal{O}.$ Therefore, for each
 an unique $f_h\in\mathcal{O}$ such that $h^{-1}[\triangle_X^f(h)]s=f_hs.$ In other words, $h^{-1}[\triangle_X^f(h)]$
 is zero order operator. This motivate the following definition.
 \begin{de}\label{2.5}{\sf A (r)-order differential operator $\varphi$ is logarithmic
 along  $D$ if $s\mapsto[\varphi(hs)-h\varphi (s)]h^{-1}$ is an
 (r-1)-order differential operators on
 $\mathcal{E}.$}
 \end{de}
 \begin{no} We denote
{\sf  $Diff^r(\mathcal{E})$  the set of $r$-order  differential
operators and $Diff^r_{\log}(\mathcal{E})$ is the subset of r-order
differential operators logarithmic along $D$ }
 \end{no}
 According to what precedes, $\triangle_X\in Diff^1_{log}(\mathcal{E});$ For all $X\in Der_X(\log D).$
 \begin{lem}\label{lem 2.7}
 Let $\varphi$ be a first order differential operator logarithmic along
 $D$ for all sections $f$ of $\mathcal{O}_X,$
There exists unique $\tilde{f}\in \mathcal{O}_X$ such that
$[\varphi(fs)-f\varphi (s)]=m_{\tilde{f}}s.$
\end{lem}
 \begin{proof}
$[s\mapsto\varphi(fs)-f\varphi (s)]\in Diff^0_{\log}(\mathcal{E})$ then
there exist $\tilde{f}\in \mathcal{O}_X$ such that
$[\varphi(fs)-f\varphi (s)]=m_{\tilde{f}}s.$ If $g$ is another
section of $\mathcal{O}_X$ such that $[\varphi(fs)-f\varphi
(s)]=m_gs,$ Then $\tilde{f}s=gs$ for all $s\in \mathcal{E}$; i.e;
$\tilde{f}=g$
\end{proof}
\begin{cor}If $\varphi$ is a first order operator logarithmic along $D$
then $\tilde{h}\in h\mathcal{O}_X$
\end{cor}
\begin{proof}
For all $s\in\mathcal{E}, \varphi(hs)-h\varphi (s)=\tilde{h}s$ and
there exist $g\in \mathcal{O}_X$ such that $\varphi(hs)-h\varphi
(s)=hgs$. Therefore, $(\tilde{h}-hg)s=0$ for all $s.$
\end{proof}
It follows that any first order differential operator logarithmic
along $D,$ $\varphi$ gives
 rise to a map
 $\sigma_{\varphi}:\mathcal{O}_X\rightarrow\mathcal{O}_X$ defined by $\sigma_{\varphi}(f)=\tilde{f}$ such that
 $[\varphi(fs)-f\varphi (s)]=\tilde{f}s$ for all $s\in\mathcal{E}.$
\begin{lem}\label{lem 2.9}
For all $\varphi\in Diff^1_{\log}(\mathcal{E}),\sigma_{\varphi}\in
H^0(X, Der^1_X(\log D))$
\end{lem}
\begin{proof}
$$\begin{array}{lll}
  \sigma_{\varphi}(f.g)s & = & \varphi(f(gs)-fg\varphi (s) \\
   & = & \sigma_{\varphi}(f)(gs)+f\varphi(gs)-fg\varphi(s) \\
   & = & \sigma_\varphi (f)(gs)+f(\varphi(gs)-g\varphi (s)) \\
   & = & (\sigma_\varphi(f)g+f\sigma_\varphi(g))s
\end{array}$$
in other hand, we have; $$\begin{array}{ccc}
                           \sigma_\varphi(h)s & = &  \varphi(hs)-h\varphi (s)\\
                           & = & hm_{\tilde{h}}(s)
                         \end{array}$$
                         Then
$(\sigma_{\varphi}(h)-hm_{\tilde{h}})s=0$ for all $s$ \\
Therefore, $\sigma_{\varphi}(h)\in h\mathcal{O}_X$; i.e.,
$\sigma_{\varphi}\in H^0(X,Der^1_X(-\log D)).$

\end{proof}

\begin{prop}\label{Pro 2.10}
$ Diff^1_{\log}(\mathcal{E})$ is closed under commutator.
\end{prop}
\begin{proof}
Let $\varphi_1, \varphi_2$ be two sections of
$Diff^1_{\log}(\mathcal{E});$
 We have:
$$\begin{array}{lll}
  \varphi_1\varphi_2(fs) & = &\varphi_1\left( f\varphi_2(s)+\bar{f}^2s\right) \\
&=& f\varphi_{1}(f\varphi_{2}(s)+\varphi_1(\bar{f}^2s))\\
&=&f\varphi_1\left( \varphi_2(s)\right) +{\bar{f}}^1\varphi_2(s)+\bar{f}^2\varphi_1(s)+{\bar{{\bar{f}}}^2}^1s
\end{array}$$
In the same way, we obtain:\\ $
\varphi_2\varphi_1(fs)=f\varphi_2\left( \varphi_1(s)\right) +{\bar{f}}^2\varphi_1(s)+\bar{f}^1\varphi_2(s)+{\bar{{\bar{f}}}^1}^2s$\\
therefore\\
$\varphi_1\varphi_2(fs)-\varphi_2\varphi_1(fs)-f\left( \varphi_1\varphi_2-\varphi_2\varphi_1\right)(s)=({\bar{{\bar{f}}}^2}^1-{\bar{{\bar{f}}}^1}^2)s $\\
In other hand, since $\varphi_1,\varphi_2$ are sections of $Diff^1_{\log}(\mathcal{E}),$ there exist $h_1, h_2\in\mathcal{O}_X$ such that:\\
$[\varphi_2(hs)-h\varphi_2(s)]\frac{1}{h}=h_2s$ and $[\varphi_1(hs)-h\varphi_1(s)]\frac{1}{h}=h_1s$\\
i.e; $\bar{h}^2=hh_2$ and $\bar{h}^1=hh_1.$\\
In the same way, there exist $h_{21},h_{12}\in\mathcal{O}_X$ such that ${\bar{\bar{h}}^1}^2=hh_{12}$ and ${\bar{\bar{h}}^2}^1=hh_{21}.$\\
Therefore,\\
$\varphi_1\varphi_2(hs)-\varphi_2\varphi_1(hs)-h\left( \varphi_1\varphi_2-\varphi_2\varphi_1\right)(s)=({\bar{{\bar{h}}}^2}^1-{\bar{{\bar{h}}}^1}^2)s=h[h_{21}-h_{12}]s$
\end{proof}
\begin{prop}\label{Pro 2.11}
$ Diff^1_{\log}(\mathcal{E})$ is a logarithmic Lie-Rinehart algebra.
\end{prop}
\begin{proof} According to above results, we have the following map.
  $$\begin{array}{ccc}
      Diff^1_{\log}(\mathcal{E}) & \rightarrow & Der_X(\log D) \\
      \varphi & \mapsto & \sigma_{\varphi}
    \end{array}$$
For all $f\in\mathcal{O}_X, s\in\mathcal{E},$ we have:

  $ \begin{array}{lllll}
      \sigma_{[\varphi_1,\varphi_2]}(f)s & = & [\varphi_1,\varphi_2](fs)-f[\varphi_1,\varphi_2](s) \\
       & = & \varphi_1\varphi_2(fs)-\varphi_2\varphi_1(fs)-f\varphi_1\varphi_2(s)+\varphi_2\varphi_2(s) \\
       & = & \varphi_1(\sigma_{\varphi_2}(f)s+f\varphi_2(s))- \varphi_2(\sigma_{\varphi_1}(f)s+f\varphi_1(s))-f[\varphi_1,\varphi_2]s\\
      & =& \varphi_1(\sigma_{\varphi_2}(f)s) + \varphi_1(f\varphi_2(s))-\varphi_2(\sigma_{\varphi_1}(f)s) - \varphi_2(f\varphi_1(s))-f[\varphi_1,\varphi_2]s\\
       & = & \sigma_{\varphi_1}(\sigma_{\varphi_2}(f))s+\sigma_{\varphi_2}(f)\varphi_1(s)+
       \sigma_{\varphi_1}(f)\varphi_2(s)+f\varphi_1(\varphi_2(s))-\\&&
      \sigma_{\varphi_2}(\sigma_{\varphi_1}(f))s-\sigma_{\varphi_1}(f)\varphi_2(s)-
       \sigma_{\varphi_2}(f)\varphi_1(s)-f\varphi_2(\varphi_1(s))- f[\varphi_1,\varphi_2]s\\
       & = & [\sigma_{\varphi_1},\sigma_{\varphi_2}](f)s
    \end{array}$
In other hand, for all $\varphi_1,\varphi_2\in Diff^1_{\log}(\mathcal{E}), f\in\mathcal{O}$ and $s\in\mathcal{E}$ we have:

$\begin{array}{lll}[\varphi_1,f\varphi_2]&=&\varphi_1(f\varphi_2(s))-(f\varphi_2)(\varphi_1(s))\\
                                         &=&f\varphi_1(\varphi_2(s))+\sigma_{\varphi_1}(f)(\varphi_2(s))-f\varphi_2(\varphi_1(s))\\
                                          &=& \sigma_{\varphi_1}(f)(\varphi_2(s))+f[\varphi_1,\varphi_2]
 \end{array}
$
\end{proof}
From above results, we deduce the  following
exact sequence of Lie-Rinehart algebras
$$\xymatrix{0\ar[r]&\mathcal{O}_X\ar[r]^m&Diff^1_{\log}(\mathcal{E})\ar[r]^\sigma&Der_X(\log D)\ar[r]&0}$$
Therefore, if we replace in (\ref{Eq7}) $\mathcal{A}$ and  $k$, respectively by $\mathcal{O}_X$ and $\mathbb{C}$ then
 Dirac
Proceeding of prequantization comes back to determine a locally free rank 1 $\mathcal{O}_X$-module $\mathcal{E}$
endowed with connection $\nabla$ such that $Diff^1_{\log}(\mathcal{E})$ is
faithful representation of $(\mathcal{O}_X,\omega).$ \\
In the following paragraph, we study the obstruction to existence of solution to this problem.

\section{Prequantization.}

  Suppose that the logsymplectic manifold $(X,\omega,
 D)$ admit a prequantum representation $(Diff^1_{\log}(\mathcal{E}),\nabla, K_{\nabla})$.\\
 For all $f,g\in H^0(X, \mathcal{O}_X)$ and $s\in \mathcal{E},$
 \[\begin{array}{llll}
   \varphi(f)\varphi(g)s & = & \varphi(f)(\varphi(g)s) \\
    & = & \varphi(f)[\nabla_{v(g)}s+2\pi igs] \\
    & = & \nabla_{v(f)}(\nabla_{v(g)}s+2\pi igs)+2\pi i(f\nabla_{v(g)}s+2\pi i fgs) \\
    & = & \nabla_{v(f)}\nabla_{v(g)}s +2\pi i\nabla_{v(f)}(gs) + 2\pi i \nabla_{v(g)}s -4\pi^2fgs \\
    & = & \nabla_{v(f)}\nabla_{v(g)}s +2\pi i (H(df).g)s + 2\pi
    ig\nabla_{v(f)}s +2\pi if\nabla_{v(g)}s -4\pi^2fg
 \end{array}\]
Changing the role of $f$ and $g$, we obtain:
\[\varphi(g)\varphi(f)s=\nabla_{v(g)}\nabla_{v(f)}s +2\pi i (H(dg).f)s + 2\pi
    ig\nabla_{v(g)}s +2\pi ig\nabla_{v(f)}s -4\pi^2gfs\]
    Therefore
    \[[\varphi(f),\varphi(g)]s=[\nabla_{v(f)},\nabla_{v(g)}]s+ 4\pi i\omega(v(f),v(g))s\]
In other hand,
\[\begin{array}{lll}
    \varphi(\{f,g\}) & = & \nabla_{v(\{f,g\})}s+2\pi i\{f,g\}s \\
     & = & \nabla_{[v(f),v(g)]}s + 2\pi i\{f,g\}s \\
     & = & [\nabla_{v(f)},\nabla_{v(g)}]-K_{\nabla}(v(f),v(g))s + 2\pi i\{f,g\}s \\
     & = & [\varphi(f),\varphi(g)]s +2\pi i\{f,g\}s - K_{\nabla}(v(f),v(g))s
  \end{array}
\]
It follow that $\varphi$ is prequantum map of $(X,\omega, D)$ if and
only if

\begin{equation}\label{E4}
    K_{\nabla}=2\pi i\omega
\end{equation}

 Since
$\omega\in H^0(X,\Omega^2_X(\log D)),$ relation \ref{E4} imply that $K_\nabla$ and
then $\nabla$ are logarithmic forms.

\begin{de}
{\sf We refer to prequantum sheaf on $(X,\omega,D)$  a rank 1 connection
$(\mathcal{M}, \nabla)$ satisfy (\ref{E4}) }
\end{de}
\subsection{Extension of prequantum sheaf.}

Our main objective being to determine existence condition of prequantum sheaf
$(\mathcal{M},\nabla)$ on $(X,D,\omega)$ satisfy (\ref{E4}), we intend in a first times to determine in which case
integral condition of $\omega$ on $X-D$ could be extended to entire $X.$
Of cause we shall know in how and when it is possible to prolong connection on $X-D$
to logarithmic connection on $D.$
First about, we recall the following of S. Litaka proved in \cite{S.I}.

\begin{prop}$\label{Prop 3}$\cite{S.I}\end{prop}{\sf
Let $F$ be a closed subset of a nonsingular variety $X$ with $F\neq
X$
 If $\omega_1$ and $\omega_2$ are rational $q-$forms such that
$\omega_1|_{X-F}=\omega_2|_{X-F}$, then $\omega_1=\omega_2$}

\begin{prop}\label{Prop 4}
 {\sf If $(\tilde{\mathcal N},\tilde{\nabla})$ is extension of a prequantum sheaf $(\mathcal N,\nabla)$ on $X^*$
$(\mathcal N,\nabla)$ and $D$ a closed reduced divisor of $X$,
then $(\tilde{\mathcal N},\tilde{\nabla})$ is a prequantum sheaf of
$X$}
\end{prop}
\begin{proof}
Since $D$ is a simple normal crossing divisor and $(\tilde{\mathcal
N},\tilde{\nabla})$ is an extension of $(\mathcal N,\nabla)$ then
$K_{\tilde{\nabla}}\mid_{X^*}=K_{\nabla}=2\pi i\omega.$ The result
follows from the proposition \ref{Prop 3}
\end{proof}
\begin{cor}\label{C1}
If $(\tilde{\mathcal N},\tilde{\nabla})$ is extension of prequantum
sheaf of $X^*$ $(\mathcal N,\nabla)$ and $D$ is simple normal
crossing, then there exists strictly close logarithmic form $\tau$
such that $\tilde{\sigma}=\sigma+\tau$
\end{cor}
\begin{proof}
From Proposition \ref{Prop 4}, we have; $d(\sigma-\sigma^0)=0.$\\
The existence of $\tau$ follow from lemma \ref{L1}.
\end{proof}
\begin{lem}
Let $(\mathcal{N},\nabla^0,K_{\nabla^0})$ be a sheaf of locally free
$\mathcal{O}_{X^*}$-module of rank $1.$ If
$(\mathcal{M},\nabla,K_{\nabla})$ be a sheaf of locally free
$\mathcal{O}_{X}$-module of rank $1$ such that\\
$\nabla=\nabla^0+\tau\otimes id$, with $\tau$ a close logarithmic
1-form, then $(\mathcal{M},\nabla,K_{\nabla})$ is prequantum sheaf
of $(X,\omega, D)$ if and only if $$K_{\nabla^0}=2\pi i\omega$$
\end{lem}
\begin{proof}
The prequantum map \ref{E4} become
\begin{equation}\label{E8}
   \begin{array}{llrr} \hat{}:\mathcal{O}_X\rightarrow End_{\mathbb{C}}(\mathcal{E})\\
    \hat{f}s=v(f).\tau\otimes s+\nabla^0 s+2\pi ifs
    \end{array}
\end{equation}
And we have by simple calculation
$$\begin{array}{llll}\hat{f}(\hat{g}s)&=v(g)\tau v(f)\tau s+
v(f)\tau\nabla^0_{v(g)}s + v(g)v(f)\tau s+ 2\pi
igv(f)\tau s + v(g)\tau\nabla^0_{v(f)}s +\\
&\nu_{v(g)}\nabla^0_{v(f)}s + 2i\pi g\nabla^0_{v(f)}s +2i\pi
v(g)\tau fs + 2i\pi f\nabla^0_{v(g)}s + 2i\pi v(g)fs-
4\pi^2gfs\end{array}$$ and
$$\begin{array}{llrr}\hat{g}(\hat{f}s)&=v(f)\tau v(g)\tau s+
v(g)\tau\nabla^0_{v(f)}s + v(f)v(g)\tau s+ 2\pi
ifv(g)\tau s + v(f)\tau\nabla^0_{v(g)}s +\\
&\nabla^0_{v(f)}\nabla^0_{v(g)}s + 2i\pi f\nabla^0_{v(g)}s +2i\pi
v(f)\tau gs + 2i\pi f\nabla^0_{v(f)}s + 2i\pi v(f)gs-
4\pi^2fgs\end{array}$$ then
$$[\hat{f},\hat{g}]s=[\nabla^0_{v(f)},\nabla^0_{v(g)}]s +
[v(f),v(g)]\tau s + 4i\pi\omega(v(f),v(g))$$ in other hand,
$$\hat{\{f,g\}}=[\nabla^0_{v(f)},\nabla^0_{v(g)}]s-K_{\nabla^0}(v(f),v(g))s+
[v(f),v(g)]\tau s + 2i\pi \omega(v(f),v(g))$$
\end{proof}
Now, we can give the main theorem of this section.
\begin{thm}\label{Th1}
Let $(X,D,\omega)$ be a log symplectic manifold such that:
\begin{enumerate}
    \item $D$ is closed reduced divisor of $X$
    \item $X-D$ is even dimensional complex sub manifold of $X$
    \item integral of $\omega$ in all closed connected surface of $X-D$ is integers multiple of $2\pi i$
    Then the the symplectic manifold $(X-D, \omega)$ is prequantizable and if
    its prequantum connection spreads on $X,$ the $(X,D,\omega )$ is prequantizable.
\end{enumerate}
\end{thm}
\begin{proof}
Since $(X-D,\omega)$ is symplectic manifold, condition 3. of Theorem is B. Kostant quantization condition.
If the prequantique sheaf of $(X-D, \omega)$ spreads on $X,$ their associated curvature coincided on $X-D$. An then
it follow from (\ref{Prop 3}) that their are equal.
\end{proof}
It follow from above result that prequantization of $X-D$ informe
us on the on of $X$.
Since $(X-D,\omega)$ is symplectic, it is prequantizable if and
only if the cohomology class $[\omega]\in
H^2(X-D,\mathbb{\mathbb{C}})$ live in $i_*(H^2(X-D,\mathbb{Z}))$.
 We shall be careful on the fat that
Obstruction coming from logarithmic De Rham cohomology and not from
De Rham cohomology of $X-D$. In general the two cohomology are not
equal. We need the Logarithmic Comparison Theorem before use only
cohomology of $X-D.$ Nowadays, it is prove that if the divisor $D$
is locally quasi-homogeny and free, the cohomology of logarithmic De
Rham complex is equal to the De Rham cohomology of the complement of
divisor. If we suppose that $D$ is closed and locally quasi-homogeny, then the $X-D$ is prequantization if
the  cohomology class of $\omega$
in $X-D$ is integer.
We can deduce the following proposition.
\begin{cor}
If $D$ is close, locally quasi-homogeny and free, if the de Rham cohomology class of $\omega$ on $X-D$ is integral,
then there exist a prequantum sheaf on $X-D.$ Besides,
 if prequantum connection on $X-D$ extending on $X,$ then $(X,D,\omega)$ heve prequantum sheaf.
\end{cor}
We remark that the problem of extending connection is fundamental in our approach. In the following paragraph, we will
study the case where $D$ is a Normal Crossing Divisor.
\subsection{The normal Crossing Divisor Case.}
Thoughout this section $X$ denotes connected complex analytic compact manifold of dimension 2n and
$D=\underset{i=1}{\overset{s}\sum}v_iD_i$
an effective normal crossing divisor on X, i.e.
an effective divisor locally with nonsingular components
meeting transversally; $\omega$ is a logsymplectic structure on $X.$
Using main result of P.Deligne and B.Malgrange about extension of connection on $X-D$
we prove the sufficient condition of prequantization of $(X,D,\omega)$. \\
First about, we recall the notion of extension of connection.
\begin{de}
If $(\mathcal{M},\nabla)$ is connection on $X-D,$ we called meromorphic extension of $\mathcal{M}$ and
$\mathcal{O}_X[D]$-coherent module $\tilde{\mathcal{M}}$ provided with isomorphism
$\tilde{\mathcal{M}}|_{X-D}=\mathcal{M}$
\end{de}
Thanks to Hilbert Nullstellensatz Theorem, B. Malgrange prove in \cite{BM} the following Lemma.
\begin{lem}\cite{BM}
A coherent $\mathcal{O}_X[D]$-module $M$ whose support is contained in $D $ is trivial.
\end{lem}
It follow that extension of connection $(\mathcal{M},\nabla)$ when it exist it is unique. Therefore,
If $(X-D, \omega)$ is symplectic manifold, then the prequantum connection of $(X,D,\omega).$
The following P.Deligne Theorem assure the existence of  extension of each connection on $X-D$ when $D$ is normal crossing divisor on X.
\begin{thm}\cite{BM}
Let $D$ be a divisor with normal crossing of $X,$ and $(\mathcal{M},\nabla)$ aconnection on $X-D.$
There exists a free extension $(G, \triangle)$ of $(\mathcal{M},\nabla)$ on $X$,
unique up unique isomorphism such that
\begin{enumerate}
    \item $\nabla$ has logarithmic pole with respect to $G$
    \item The eigenvalues of the residues of $\nabla$ with respect to $G$ belong to the image of
$\tau$ the section of $\mathbb{C}\rightarrow \mathbb{C}/\mathbb{Z}$

\end{enumerate}

\end{thm}
We can then state the main theorem of this section.
\begin{thm}\label{Th2}
Let $(X,D,\omega)$ be a log symplectic manifold such that:
\begin{enumerate}
    \item $D$ is normal crossing divisor of $X$
    \item $X-D$ is even dimensional complex sub manifold of $X$
    \item integral of $\omega$ in all closed connected surface of $X-D$ is integers multiple of $2\pi i$
    Then  $(X,D,\omega )$ is prequantizable.
\end{enumerate}

\end{thm}
 when it exist is unique. According to above results, if extension of each connection on $X-D$ exist, then
existence of prequantum connection on $X-D$ imply that of $X.$ Since the unique extension of prequantum
connection $(\mathcal{M},\nabla)$ on $X-D$ will agree on $X-D$ with $(\mathcal{M},\nabla)$.
\begin{ex}
If $\omega$ is exact, then $(X,\omega, D)$ is prequantizable.
\end{ex}
\begin{proof}
If $\omega$ is exact, then $[\omega]\in H^2(X-D,\mathbb C)$ is
interger. Therefore the Kostant theorem imply that there exist
integrable connection $(\mathcal{N},\nabla)$ on $X^*$ and from
Deligne theorem,  $(\mathcal{N},\nabla)$ extend to
$(\tilde{\mathcal{N}},\tilde{\nabla})$ where the $\tilde{\nabla})$
is integrable logarithmic connection on $X.$
\end{proof}

\section{Lie algebroid formalism}
\begin{de}\cite{AP}
A Lie Algebroid on $X$ is an $\mathcal{O}_X$-module $L$ equipped
with a Lie algebra bracket $[.,.]$ and an $\mathcal{O}_X$-linear
morphism of Lie algebras $\rho:L\rightarrow T_X$ such that for $l_1,
l_2\in \mathcal{O}_X$ one has
\begin{equation}\label{E9}
    [l_1, fl_2]=f[l_1,l_2]+\rho(l_1)(f)l_2
\end{equation}
$\rho$ is called anchor.
\end{de}
We remark that Lie-Rinehart algebra is Lie-algebroid on affine
scheme. \\There exists many examples of Lie-algebroid in literature:
\begin{ex}
\begin{enumerate}
    \item Let $X$ be a smooth scheme. Then sheaf $T_X$ of tangent
    vector fields is Lie- algebroid with anchor $Id_{T_X}.$
    \item Any real Lie-algebra $g$ is Lie-algebroid on
    $X=\{\bullet\}$ the associated anchor is zero map.
    \item Let $(X, D)$ be a logarithmic manifold. The sheaf $T_X(-\log
    D)$ of logarithmic vector fields endowed with the inclusion
    morphism $T_X(-\log D)\hookrightarrow T_X$

\end{enumerate}
\end{ex}
The notion of module on Lie-Rinehart algebra is generalize by
\begin{de}\cite{AP} Let $(L,\rho, [., .])$ be a Lie-algebroid
An $L$-module is an $\mathcal{O}_X$-module $M$ equipped with Lie
action; $\omega:L\rightarrow End_{\mathcal{O}_X}(M)$ which is a
homomorphism of Lie-algebra; such that for all $f\in \mathcal{O}_X,
l\in L$ and $x\in M:$
\begin{enumerate}
    \item $(\omega(fl))(x)=f((\omega(l))(x))$
    \item $(\omega(l))(fx)=\rho(l)(f)x + (\omega(fl))(x)$
\end{enumerate}
\end{de} To simplify the notation we will denoted
$(\omega(l))(x)=l(x)$ \\
\begin{de} Let $(L,\rho,[.,.])$ be a Lie-algebroid.
Universal objet of $L$ is a sheaf of $\mathcal{O}_X$-algebras $U(L)$
equipped with a morphism $i_{\mathcal{O}_X}:\mathcal{O}_X\rightarrow
U(L)$ of $\mathcal{O}_X$-algebras and a morphism $i_L:L\rightarrow
U(L)$ of Lie-algebras having the following properties

\begin{equation}\label{E10}
    [i_L(l),i_{\mathcal{O}_X}(x)]=i_{\mathcal{O}_X}(l(x))\quad\quad
i_{\mathcal{O}_X}(f)i_{L}(l)=i_L(fl)
\end{equation}
and $(U(L),i_{\mathcal{O}_X},i_L)$ is universal among triplet
$(V,\alpha, \beta)$ satisfy \ref{E10}.
\end{de}
The notion of universal enveloping algebra of Lie-algebroid is very
useful when we define Poisson cohomology in the framework of
Algebraic Geometry, we refer the reader to \cite{AP} for more
explanation.\\
Another notion that we need is extension of lie-algebroid.
\begin{de}
Let $L$ be Lie-algebroid on $X.$ An extension of $L$ by and
$\mathcal{O}_X$-module  $F$ is an exact sequence of Lie-algebroids
\begin{equation}\label{E11}
    \xymatrix{0\ar[r]&F\ar[r]^i&E\ar[r]^\pi&L\ar[r]&0}
\end{equation}
where $E$ is Abelian Lie-algebroid.
\end{de}
Many notion are related to Lie-algebroid extension. For example:\\
A transverse of an extension of Lie-algebroid $L$ is a morphism of
$\mathcal{O}_X$-modules $\chi:L\rightarrow E$ such that $\pi\circ
\chi$.\\
A back-transverse of an extension of Lie-algebroid is a morphism of
$\mathcal{O}_X$-modules $\lambda:E\rightarrow L'$ such that
$\lambda\circ i=Id_E$\\
A transverse is flat if it homomorphism of Lie-algebroid.
\begin{prop}
There is one to one correspondence between the transverses $\chi$
and back-transverses $\lambda$ giving by
\begin{equation}\label{E12}
    i\circ \lambda + \chi\circ \pi=Id_E \quad\quad such\quad \quad
    that\quad\quad \lambda\circ \chi=0
\end{equation}
\end{prop}
\begin{proof}
It is an adaptation of proof in the smooth case.
\end{proof}
When $\pi$ is a submersion one to of fiber bundles,transverse all
way exist and the choice of the transverse determine an isomorphism
$E\backsimeq F\bigoplus L$.  The following proposition is an
generalization of Proposition 2.13 (see \cite{JH})
\begin{prop}\label{Prop7}
For any Lie-algebroid $L$ and an $\mathcal{O}_X$-module locally free
$\mathcal{E},$ there is up to congruence extensions at most one
extension of $L$ by $End_{\mathcal{O}_X}(\mathcal{E})$
\end{prop}
\begin{proof}
Let $(L,[-,-]_\rho,\rho)$ be an Lie algebroide on $X$ and
$\mathcal{M}$ an $\mathcal{O}_X$-module. Define on
$\mathcal{E}nd_{\mathbb{C}}(\mathcal{M})\oplus\mathcal{L}$ the
following bracket

\begin{equation}\label{E13}
   [(\beta,l),(\beta',l')]=(\beta\beta'-\beta'\beta,[l,l']_\rho)
\end{equation}

For each open subset $U$ of $X$ define

$\Gamma(U,\mathcal{A}(\mathcal{M}))=\{(\beta,l)\in
\mathcal{E}nd_{\mathbb{C}}(\mathcal{M})\oplus\mathcal{L},
\beta(fm)=\rho(l)(f)m+f\beta(m)\}$.

\textbf{$\mathcal{A}$, is a sheaf of $\mathcal{O}_X$-module}\\
Indeed, for all $f\in \mathcal{O}_X, (\beta,
l)\in\mathcal{A}(\mathcal{M})$.
Define $f(\beta,l):=(f\beta,fl)$.\\
Since \\
\[\begin{array}{llll}
        f\beta(gm) & = & f(\beta(gm)) &  \\
         & = & f(\rho(l)(g)m+g\beta(m)) &  \\
         & = & f\rho(l)(g)m + fg\beta(m) &  \\
         & = & \rho(fl)(g)m + fg\beta(m) &  \\
         & = & \rho(fl)(g)m+g((f\beta)(m)) &
      \end{array}\]

 Then $f(\beta,l)\in \mathcal{A}(\mathcal{M}).$

\textbf{$\mathcal{A}$, is a Lie algebroid.}\\
Indeed, let $(\beta,l),(\beta',l')\in\mathcal{A}(\mathcal{M})$.\\
We have:
\[\begin{array}{llll}
  \beta\beta'(fm) & = & \beta(\beta(fm)) &  \\
   & = & \beta(\rho(l')(f)m+f\beta'(m)) &  \\
   & = & \beta(\rho(l')(f)m)+\beta(f\beta'(m)) &  \\
   & = & \beta(\rho(l')(f)m)+ \rho(l)(f)\beta'(m)+f(\beta\beta'(m)) &  \\
   & = & \rho(l)(\rho(l')(f))m + \rho(l')(f)(\beta(m))+ \rho(l)(f)\beta'(m)+f(\beta\beta'(m))&
\end{array}\]
In the same manner, we obtain\\ $\beta'\beta(fm)=
\rho(l')(\rho(l)(f))m + \rho(l)(f)(\beta'(m))+
\rho(l')(f)\beta(m)+f(\beta'\beta(m))$\\ Then
$[\beta,\beta'](fm)=(\rho[l,l']_\rho(f))m+f([\beta,\beta'](m)).$\\
Therefore,
$[(\beta,l),(\beta',l')]=(\beta\beta'-\beta'\beta,[l,l']_\rho)\in
\mathcal{A}(\mathcal{M})$

\textbf{Leibniz property}\\
Let $(\beta,l),(\beta',l')\in \mathcal{A}(\mathcal{M}),
f\in\mathcal{O}_X$. $[(\beta,l),f(\beta',l')]=(\beta
f\beta'-f\beta'\beta,[l,fl']_\rho).$\\ We have
$\beta(f\beta')(m)-(f\beta')\beta(m)=\rho(l)(f)\beta'(m)+f(\beta\beta'(m)-\beta'\beta(m)).$\\
Then
$\beta(f\beta')-(f\beta')\beta=\rho(l)(f)\beta'+f[\beta,\beta'].$\\
Then $(\beta
f\beta'-f\beta'\beta,[l,fl']_\rho)=\rho(l)(f)(\beta',l')+f()[\beta,\beta'],[l,l']_\rho.$
\\ Therefore the Lie Leibniz property is satisfy we anchor $\Phi(\beta,l)=\rho(l).$
\end{proof}
\begin{prop}\label{Prop8}
\begin{enumerate}
    \item For all $(\beta,l)\in\mathcal{A}(\mathcal{M})$, $\beta\in Diff^1(\mathcal{M})$
    \item If $L=Der(-\log D)$ then $\beta\in Diff^1_{log}(\mathcal{M})$
    \item the map $\varphi((\beta,l))=\beta$ is Lie-algebroid
    homomorphism
\end{enumerate}

\end{prop}
\begin{proof}
\begin{enumerate}
    \item Let $(\beta,l)\in\mathcal{A}(\mathcal{M}).$
    $\beta(fs)=(\rho(l)(f))s+f(\beta(s))$\\ Then
    $\beta(fs)-f(\beta(s))=(\rho(l)(f))s$\\ Therefore $\beta\in Diff^1(\mathcal{M})$
    \item Suppose that $L=Der(-log D)$, then $\rho(l)\in Der(-log
    D)$ and then $\rho(l)(h)\in h\mathcal{O}_X$ for all $l\in Der(-log
    D).$\\ Therefore $[\beta(fs)-f(\beta(s))]h^{-1}\in\mathcal{O}_X$
    \item consider the map
    $$\begin{array}{cccc}
      \varphi: & \mathcal{A}(\mathcal{M}) & \rightarrow & Diff^1(\mathcal{M}) \\
       & (\beta,l) & \mapsto & \beta
    \end{array}$$
    $\varphi$ is an homomorphism of Lie-algebroids.\\
    Indeed for all $(\beta,l)\in\mathcal{A}(\mathcal{M}), s\in \mathcal{M};$
    $$\begin{array}{lll}
       \sigma\circ\varphi((\beta,l))(f)(s) & = & \sigma_{\varphi(\beta,l)}(f)(s) \\
        & = & \sigma_\beta(f)(s) \\
        & = & \beta(fs)-f\beta(s) \\
        & = & (\rho(l)(f))(s)\\
        & = & (\Phi(\beta,l)(f))(s)
     \end{array}
    $$
    then $\sigma\circ\varphi=\Phi$ where $\Phi(\beta,l)=\rho(l)$
\end{enumerate}
\end{proof}
\begin{cor}
For all $\mathcal{O}_X$-module locally free of rang 1,
$\mathcal{M},$
$Diff^1_{\log}(\mathcal{M}),\mathcal{A}_{\log}(\mathcal{M})\in
Ext^1(\mathcal{O}_X,Der(-\log D))$
\end{cor}
\begin{proof}
It follow from the five lemma and above corollary.
\end{proof}
It follows from this proposition that giving a complex line bundle
$L$ on a logarithmic manifold $(X,D),$ There exist up to congruence
of extensions at most one extension of the Lie algebroid $T_X(-\log
D)$
\begin{cor}For all $\mathcal{O}_X$-module locally free of rank one
$L$ there exist an exact sequence of Lie algebras
\begin{equation}\label{E13}
    \xymatrix{0\ar[r]&End(\Gamma(\mathcal{E}))\ar[r]^i&A_{\log}(\mathcal{E})\ar[r]^\pi&Der(-\log D)\ar[r]&0}
\end{equation}
  \end{cor}
  Therefore, we have the following definition
  \begin{de}
An extension $A_{\log}(\mathcal{E})$ of $Der(-\log D)$ giving by
relation  (\ref{E13}) is called Atiyah logarithmic algebroid of the
invertible sheaf $\mathcal{E}.$
  \end{de}
The existence of $A_{\log}(\mathcal{E})$ allowed us to think about
representation of logsymplectic  Poisson algebra $(\mathcal{O}_X,
\omega)$ by $A_{\log}(\mathcal{E}).$ \\
In this case prequantum representation shall commuted the following
diagram.
\begin{equation}\label{D1}
   \xymatrix{0\ar[r]&\mathcal{O}_X\ar[r]&A_{\log}(\mathcal{E})\ar[r]&Der(-\log D)\ar[r]&0\\
         0\ar[r]&\mathbb{C}\ar[r]\ar[u]&(\mathcal{O}_X,\omega)\ar[u]^Q\ar[r]&Ham(X)\ar[u]\ar[r]&0}
\end{equation}
\begin{lem}Let $\nabla:Der(\log D)\rightarrow End(\mathcal{M})$ be a
logarithmic connection on $\mathcal{M}$. For all $\delta\in
Der(-\log D)$
\begin{enumerate}
    \item $\sigma_{\nabla_\delta}=\delta$
    \item For all $\varphi\in \mathcal{A}_{\log}(\mathcal{M}), \nabla_{\sigma_\varphi}-\varphi\in \ker(\sigma)\backsimeq\mathcal{O}_X$
\end{enumerate}
\end{lem}
\begin{proof}
Let $\nabla$ be a logarithmic connection on $\mathcal{M}$ and
$\delta\in Der(-\log D).$
\begin{enumerate}
    \item $\nabla_{\delta}(fs)=f\nabla_{\delta}s+ \delta(f)s$; i.e.
    $\delta(f)s=\nabla_{\delta}(fs)-f\nabla_{\delta}s$.\\ Therefore,
    there exist $g\in \mathcal{O}_X$ such that
    $[\nabla_{\delta}(hs)-h\nabla_{\delta}s]h^{-1}=gs;$ i.e. $\nabla_\delta\in \mathcal{A}_{\log D}(\mathcal{M})$
    we can compute its image by $\sigma.$ We have:
    $$=\begin{array}{llll}
                                   \sigma_{\nabla_\delta}(f)(s) & = & \nabla_\delta(fs)-f\nabla_\delta s \\
                                     & = & (\delta(f))s \quad\forall f\in\mathcal{O}_X, s\in \mathcal{M} \\

                                  \end{array}
    $$
    Then $\sigma_{\nabla_\delta}=\delta$
    \item $\forall\varphi\in \mathcal{A}_{\log}(\mathcal{M})$,
    $$\begin{array}{lll}
      \sigma_{\nabla_{\sigma_\varphi}}(f)(s) & = & \nabla_{\sigma_\varphi}(fs)-f\nabla_{\sigma_\varphi}s\\
       & = & f\nabla_{\sigma_\varphi}s+(\sigma_\varphi(f))s-f\nabla_{\sigma_\varphi}s \\
       & = &  (\sigma_\varphi(f))s\\
      i.e.,\,\,  \sigma_{\nabla_{\sigma_\varphi}}(f)(s)& = & \sigma_\varphi(f)s\forall f\in\mathcal{O}_X, s\in\mathcal{M} \\
     i.e.,\,\,\,  \sigma_{\nabla_{\sigma_\varphi}}- \sigma_\varphi & = & 0 \\
     i.e.,\,\,\,  \sigma_{\nabla_{\sigma_\varphi}- \varphi}  & = & 0 \\
      i.e.,\,\,\, \nabla_{\sigma_\varphi}- \varphi & \in  & \ker(\sigma) \\

    \end{array}$$
\end{enumerate}
\end{proof}
It follow that there exist $m(\varphi)\in\mathcal{O}_X$ such that
\begin{equation}\label{E14}
\varphi=\nabla_{\sigma_\varphi}+m(\varphi)\end{equation}

\begin{lem}
Let $\omega$ be a logsymplectic form on $X,$ the following morphism
$Der(-\log D)\rightarrow \Omega^1_X(\log D)\quad by\quad
\delta\mapsto i_{\delta}\omega$ is an isomorphism if $D$ is free.
\end{lem}
It follow from this lemma that $\forall f\in\mathcal{O}_X$ there
exist $\delta_f\in Der_X(-\log D)$ such that $i_{\delta_f}\omega=df$
\begin{lem}There exist a map $\lambda:(\mathcal{O}_X,\omega)\rightarrow \mathcal{A}_{\log}(\mathcal{M})$
who commute the following diagram
$$\xymatrix{(\mathcal{O}_X,\omega)\ar[d]_\gamma\ar[r]^\lambda&\mathcal{A}_{\log}(\mathcal{M})\ar[dl]^\sigma\\Der(-\log
D)}$$
\end{lem}
\begin{proof}
For all $f\in\mathcal{O}_X$ since $\sigma$ is onto, it admit a
section $\tau$ such that $\sigma\circ\tau=Id$ we denoted
$\lambda=\tau\circ\gamma$
\end{proof}
For all $f\in \mathcal{O}_X, \lambda_f$ satisfy equation (\ref{E14})
i.e. $\lambda_f=\nabla_{\sigma_{\lambda_f}}+m(f)$ where we have
replaced $m(\lambda_f)$ by $m(f).$
\begin{cor}$m$ is $\mathcal{O}_X$-linear on $\mathcal{M}$ \end{cor}
\begin{proof}We have $$\begin{array}{lll}
                        \lambda_f(gs) & = & \nabla_{\sigma_{\lambda_f}}(gs)+m(f)(gs) \\
                         & = & g\nabla_{\sigma_{\lambda_f}}s+\sigma_{\lambda_f}(g)s+m(f)(gs) \\
                        \lambda_f(gs)-g\lambda_f(s) & = & g\nabla_{\sigma_{\lambda_f}}s+\sigma_{\lambda_f}(g)s+m(f)(gs)-g\nabla_{\sigma_{\lambda_f}}s+g\sigma_{\lambda_f}  \\
                        \lambda_f(gs)-g\lambda_f(s) & = &
                        \sigma_{\lambda_f}(g)(s)+
                        (m(f))(gs)-g(m(f))(s)\\
 \i.e. (m(f))(gs)-g(m(f))(s)&=&0
                      \end{array}
$$
\end{proof}

\begin{cor}\begin{enumerate}
    \item The map $\mathcal{Q}$ in diagram (\ref{D1}) $\mathcal{Q}:(\mathcal{O}_X,\omega)\rightarrow \mathcal{A}_{\log}(\mathcal{M})$
have the form $\mathcal{Q}(f)=\nabla_{\delta_f}+\alpha m(f)$ for
some constant $\alpha$
    \item $\mathcal{Q}(f)=\nabla_{\delta_f}+\alpha m(f)$ is
    prequantization representation if and only if
    $m:(\mathcal{O}_X,\omega)\rightarrow(\mathcal{O}_X,\omega)$ is
    $\mathbb{C}$-linear and for any $f,g\in (\mathcal{O}_X,\omega)$
and
\begin{equation}\label{E15} \dfrac{K_\nabla(\delta_f,\delta_g)}{\alpha}=\delta_gm(f)-\delta_fm(g)+m\{f,g\}\end{equation}
\end{enumerate}
\end{cor}
Let us recall that the Poisson structure $\{-,-\}_{\omega}$ induced by $\omega$ induce on $\mathcal{O}_X$ a structure of $\mathcal{O}_X$-module defining by $f.g:=\delta_f.g=\{f,g\}.$ So, we can consider the cohomology of the complex $\mathcal{L}alt_{\log}^*(\mathcal{O}_X)=\underset{k=0}{\overset{\infty}\sum}
\mathcal{L}alt_{\log}^k(\mathcal{O}_X)$ where $\mathcal{L}alt_{\log}^k(\mathcal{O}_X)=\{m:\mathcal{O}_X\times...\times\mathcal{O}_X\rightarrow \mathcal{O}_X
\quad \texttt{alternating multilinear}\}$ the associated differential is :
$d_{\log}m(f_0,...,f_r)=\underset{k=0}{\overset{k}\sum}f_km(f_0,...,\hat{f_k},...,f_r)+ \underset{k\leq j}{\overset{r}\sum}m(\{f_k,f_j\},f_0,...,\hat{f_k},...,\hat{f_j},...,f_r)$ and $fm(g):=\{f,m(g)\}$
Moreover, any logarithmic $r$-form $\eta\in\Omega^r(\log D)$ define a r-cochain
 $\mathcal{K}_{\eta}\in\mathcal{L}alt_{\log}^r(\mathcal{O}_X)$ via $\mathcal{K}_{\eta}(f_1,...,f_r):=\eta(\delta_{f_1},...,\delta_{f_r})$ \\ It follow that $m$ is 1-cochaine and the condition that $\mathcal{Q}$ be a representation then becomes

\begin{equation}\label{E16}
d_{\log}m(f,g)=\alpha\mathcal{K}_{K_\nabla}(f,g)
\end{equation}

\begin{prop}The map $\mathcal{K}:\Omega^*(\log D)\rightarrow\mathcal{L}alt_{\log}^*(\mathcal{O}_X)$ $\eta\mapsto\mathcal{K}_\eta$ preserves the wedge product and is an injection of cochain complex when $D$ is free
\end{prop}
before giving the integrally condition of $\omega$, we shall first recall the following useful notions
\begin{de}
\begin{enumerate}
 \item The polynomial $h(z_1,...,z_n)=\sum a_{i_1...i_n}z^{i_1}...z^{i_n}\in\mathcal{O}_{C^n}$ is weighted homogeneous if there exist positive integer weights $w_1,...,w_n$ such that $h(z_1^{w_1},...,z_n^{w_n})$ is homogeneous.
\item The divisor $D\subset X$ is locally quasi-homogeneous if for all $x\in D$ there are local coordinates on $X$, centered at $x$, with respect to which $D$ has a weighted homogeneous defining equation.
\end{enumerate}
\end{de}
\begin{prop}\cite{NA}\label{Prop9}
Let $D$ be a strongly quasihomogeneous free divisor in the complex
manifold $X$, let $U$ be the complement of $D$ in $X,$ and let
$j:U\rightarrow X$ be inclusion. Then the natural morphism from the
complex $\Omega^*_X(\log D)$ of differential forms with logarithmic
poles along $D$ to $Rj*\mathbb{C}$ is quasi-isomorphism.
\end{prop}
\begin{ex}\end{ex}
\begin{enumerate}
    \item Let $X=\mathbb{C}^3$ and $D$ the divisor of $X$ defined by the
equation \\
$xy(x+y)((x-2)x+y)=0, $\\
$Der_X(-\log D)$ has free basis $\{\delta_1,\delta_2,\delta_3\}$\\
$\delta_1=x\partial_x+y\partial_y$\\
$\delta_2=((z-2)x+y)\partial_z$\\
$\delta_3=x^2\partial_x-y^2\partial_y-(z-2)(x+y)\partial_z$. Each
$\delta_i$ is logarithmic vector field and the determinant of their
coefficients is reduced equation of $D.$ It follow from K.Saito
theorem that the system $(\delta_1,\delta_2, \delta_2)$ is a basis
of $Der(-\log D).$ Therefore, $D$ is free. Since no linear
combination of $(\delta_1,\delta_2, \delta_2)$ has non-singular
linear part, $D$ can not be quasihomogeneous.
    \item All normal crossing divisor is locally quasihomogeneous.
\end{enumerate}
It follow from \ref{Prop9} that and Grothendieck's Comparison
Theorem that Cohomology of $X-D$ compute the one of the complex
$(\Omega_X^*(\log D),d).$ We denote $H^*_{dR\log}(X)$ the cohomology
of the complex $(\Omega_X^*(\log D),d)$\\
Let $X$ be a complex analytic space, $\mathcal{F}$ a coherent sheaf
on $X.$ Denote by $\mathfrak{G}_k(\mathcal{F}):=\{m\in X;
prof_m\mathcal{F}\leq k\}.$ \\ We saying that the sheaf

$\mathcal{F}$ satisfy the condition $(s_k)$ if
$$\dim\mathfrak{G}_k(\mathcal{F})\leq k-2$$
\begin{thm}If $D$ is zero dimensional locally homogeneous free divisor of $X$ and if the De Rham cohomology
class of $\omega$ on $X-D$ live in $i_*(H^2(X-D),\mathbb{Z})$ then
$(X,\omega)$ have prequantum bundle if the associated prequantum
bundle of $X-D$ satisfy the condition $(s_2)$
 \end{thm}
 \begin{proof}
Since the De Rham cohomology class of $\omega$ on $X-D$ live in
$i_*(H^2(X-D)$, it follow from B.Kostant in [\cite{B.K}] that there
exist a rank one locally free $\mathcal{O}_{X-D}$-module
$\mathcal{F}$ such that the curvature satisfy the equation\ref{E16}
with $\alpha=-2\pi i.$ If $\mathcal{F}$ satisfy the condition
$(s_2)$ and $D$ is zero dimensional analytic divisor of $X,$ then
according to Trautmann Theorem that $\mathcal{F}$ there exist an
unique analytic coherent sheaf $\tilde{\mathcal{F}}$ on $X$
extending $\mathcal{F}.$ Since the curvature of $\mathcal{F}$
coincide on $X-D$ with curvature of $\tilde{\mathcal{F}}$ it follow
from Proposition \ref{Prop 3} and to the Logarithmic Comparison
Theorem that $\tilde{\mathcal{F}}$ is prequantum sheaf of $X$
 \end{proof}





\label{lastpage}
\end{document}